# $\gamma$ and $\beta$ approximations via general ordered topological spaces

M. Abo-Elhamayel
Mathematics Department, Faculty of Science
Mansoura University

**Abstract**

In this paper, we introduce the concepts of $\gamma$ and $\beta$ approximations via general ordered topological approximation spaces. Also, increasing (decreasing) $\gamma$, $\beta$ boundary, positive and negative regions are given in general ordered topological approximation spaces (GOTAS, for short). Some important properties of them were investigated. From this study, we can say that studying any properties of rough set concepts via GOTAS is a generalization of Pawlak approximation spaces and general approximation spaces.

## 1. Introduction

Rough set theory was first proposed by Pawlak for dealing with vagueness and granuality in information systems. Various generalizations of Pawlak's rough set have been made by replacing equivalence relations with kinds of binary relations and many results about generalized rough set with universe being finite were obtained [4, 13,14,15,16,17,18]. An interesting and natural research topic in rough set theory is studing it via topology [5,12]. Neighborhood systems was first applied in generalizing rough sets in 1998 by T. Y. Lin as a generalization of topological connections with rough sets. Lin also introduced the concept of granular computing as a form of topological generalizations [36,37,38,39]. In this paper we give the concept of $\gamma$, $\beta$ via topological ordered spaces and studied their properties which may be viewed as a generalization of previous studies in general approximation spaces, as if we take the partially ordered relation as an equal relation we obtain the concepts in general approximation spaces[2].

## 2. Preliminaries

In this section, we give an account for the basic definitions and preliminaries to be used in the paper.

**Definition 2.1[10].** A subset $A$ of $U$, where $(U,\rho)$ is a partially ordered set is said to be increasing (resp. decreasing) if for all $a \in A$ and $x \in U$ such that $a\rho x$ (resp. $x\rho a$) imply $x \in A$.

**Definition 2.2[10].** A triple $(U,\tau,\rho)$ is said to be a topological ordered space, where $(U,\tau)$ is a topological apace and $\rho$ is a partially order relation on $U$.

**Definition 2.3[11].** Information system is a pair $(U,\mathbf{A})$, where $U$ is a non-empty finite set of objects and $\mathbf{A}$ is a non-empty finite set of attributes.

**Definition 2.4[4].** A non-empty set $U$ equipped with a general relation $R$ which generate a topology $\tau_R$ on $U$ and a partially order relation $\rho$ written as $(U, \tau_R, \rho)$ is said to be general ordered topological approximation sapce (for short, GOTAS).

**Definition 2.5[3].** Let $(U, \tau_R, \rho)$ be a GOTAS and $A \subseteq U$. We define:
  (1) $\underline{R}_{Inc}(A) = A^{\circ Inc}$, $A^{\circ Inc}$ is the greatest increasing open subset of $A$.
  (2) $\underline{R}_{Dec}(A) = A^{\circ Dec}$, $A^{\circ Dec}$ is the greatest decreasing open subset of $A$.
  (3) $\overline{R}^{Inc}(A) = \overline{A}^{Inc}$, $\overline{A}^{Inc}$ is the smallest increasing closed superset of $A$.
  (4) $\overline{R}^{Dec}(A) = \overline{A}^{Dec}$, $\overline{A}^{Dec}$ is the smallest decreasing closed superset of $A$.

  (5) $\alpha^{Inc} = \dfrac{card(\underline{R}_{Inc}(A))}{card(\overline{R}^{Inc}(A))}$ (resp. $\alpha^{Dec} = \dfrac{card(\underline{R}_{Dec}(A))}{card(\overline{R}^{Dec}(A))}$) and $\alpha^{Inc}$ (resp. $\alpha^{Dec}$)

is $R-$ increasing (resp. decreasing) accuracy.

**Definition 2.6[4].** Let $(U, \tau_R, \rho)$ be a GOTAS and $A \subseteq U$. We define:
  (1) $\underline{S}_{Inc}(A) = A \cap \overline{R}^{Inc}(\underline{R}_{Inc}(A))$, $\underline{S}_{Inc}(A)$ is called $R-$ increasing semi lower.
  (2) $\overline{S}^{Inc}(A) = A \cup \underline{R}_{Inc}(\overline{R}^{Inc}(A))$, $\overline{S}^{Inc}(A)$ is called $R-$ increasing semi upper.
  (3) $\underline{S}_{Dec}(A) = A \cap \overline{R}^{Dec}(\underline{R}_{Dec}(A))$, $\underline{S}_{Dec}(A)$ is called $R-$ decreasing semi lower.
  (4) $\overline{S}^{Dec}(A) = A \cup \underline{R}_{Dec}(\overline{R}^{Dec}(A))$, $\overline{S}^{Dec}(A)$ is called $R-$ decreasing semi upper.

$A$ is $R-$ increasing (resp. decreasing) semi exact if $\underline{S}_{Inc}(A) = \overline{S}^{Inc}(A)$ (resp. $\underline{S}_{Dec}(A) = \overline{S}^{Dec}(A)$), otherwise $A$ is $R-$ increasing (resp. decreasing) semi rough.

**Proposition 2.7[1].** Let $(U, \tau_R, \rho)$ be a GOTAS and $A \subseteq U$. Then:
  1- $\underline{R}_{Inc}(A) \subseteq \underline{\alpha}_{Inc}(A) \subseteq \underline{S}_{Inc}(A)$ $(\underline{R}_{Dec}(A) \subseteq \underline{\alpha}_{Dec}(A) \subseteq \underline{S}_{Dec}(A))$.
  2- $\overline{S}^{Inc}(A) \subseteq \overline{\alpha}^{Inc}(A) \subseteq \overline{R}^{Inc}(A)$ $(\overline{S}^{Dec}(A) \subseteq \overline{\alpha}^{Dec}(A) \subseteq \overline{R}^{Dec}(A))$.

## 3. New approximations and their properties

In this section, we introduce some definitions and propositions about near approximations, near boundary regions via GOTAS which are essential for a present study.

**Definition 3.1.** Let $(U, \tau_R, \rho)$ be a GOTAS and $A \subseteq U$. We define:

(1) $\underline{\gamma}_{Inc}(A) = A \cap [\overline{R}^{Inc}(\underline{R}_{Inc}(A)) \cup \underline{R}_{Inc}(\overline{R}^{Inc}(A))]$, $\underline{\gamma}_{Inc}(A)$ is called $R$-increasing $\gamma$ lower.

(2) $\overline{\gamma}^{Inc}(A) = A \cup [\overline{R}^{Inc}(\underline{R}_{Inc}(A)) \cup \underline{R}_{Inc}(\overline{R}^{Inc}(A))]$, $\overline{\gamma}^{Inc}(A)$ is called $R$-increasing $\gamma$ upper.

(3) $\underline{\gamma}_{Dec}(A) = A \cap [\overline{R}^{Dec}(\underline{R}_{Dec}(A)) \cup \underline{R}_{Dec}(\overline{R}^{Dec}(A))]$, $\underline{\gamma}_{Dec}(A)$ is called $R$-decreasing $\gamma$ lower.

(4) $\overline{\gamma}^{Dec}(A) = A \cup [\overline{R}^{Dec}(\underline{R}_{Dec}(A)) \cup \underline{R}_{Dec}(\overline{R}^{Dec}(A))]$, $\overline{\gamma}^{Dec}(A)$ is called $R$-decreasing $\gamma$ upper.

$A$ is $R$-increasing (resp. $R$-decreasing) $\gamma$ exact if $\underline{\gamma}_{Inc}(A) = \overline{\gamma}^{Dec}(A)$ (resp. $\underline{\gamma}_{Dec}(A) = \overline{\gamma}^{Dec}(A)$) otherwise $A$ is $R$-increasing (resp. $R$-decreasing) $\gamma$ rough.

**Proposition 3.2.** Let $(U, \tau_R, \rho)$ be a GOTAS and $A, B \subseteq U$. Then

(1) $A \subseteq B \to \overline{\gamma}^{Inc}(A) \subseteq \overline{\gamma}^{Inc}(B)$ ( $A \subseteq B \to \overline{\gamma}^{Dec}(A) \subseteq \overline{\gamma}^{Dec}(B)$ ).

(2) $\overline{\gamma}^{Inc}(A \cap B) \subseteq \overline{\gamma}^{Inc}(A) \cap \overline{\gamma}^{Inc}(B)$ ($\overline{\gamma}^{Dec}(A \cap B) \subseteq \overline{\gamma}^{Dec}(A) \cap \overline{\gamma}^{Dec}(B)$).

(3) $\overline{\gamma}^{Inc}(A \cup B) \supseteq \overline{\gamma}^{Inc}(A) \cup \overline{\gamma}^{Inc}(B)$ ($\overline{\gamma}^{Dec}(A \cup B) \supseteq \overline{\gamma}^{Dec}(A) \cup \overline{\gamma}^{Dec}(B)$).

**Proof.**

(1) Omitted.

(2) $\overline{\gamma}^{Inc}(A \cap B) = (A \cap B) \cup [\overline{R}^{Inc}(\underline{R}_{Inc}(A \cap B) \cup \underline{R}_{Inc}(\overline{R}^{Inc}(A \cap B)))]$

$\subseteq (A \cap B) \cup [\overline{R}^{Inc}(\underline{R}_{Inc}(A) \cap \underline{R}_{Inc}(B)) \cup \underline{R}_{Inc}(\overline{R}^{Inc}(A) \cap \overline{R}^{Inc}(B)))]$

$\subseteq (A \cap B) \cup [\overline{R}^{Inc}(\underline{R}_{Inc}(A) \cap \overline{R}^{Inc}\underline{R}_{Inc}(B)) \cup \underline{R}_{Inc}(\overline{R}^{Inc}(A) \cap \underline{R}_{Inc}\overline{R}^{Inc}(B)))]$

$\subseteq A \cup [\overline{R}^{Inc}(\underline{R}_{Inc}(A) \cup \underline{R}_{Inc}(\overline{R}^{Inc}(A)] \cap B \cup [\overline{R}^{Inc}(\underline{R}_{Inc}(B) \cup \underline{R}_{Inc}(\overline{R}^{Inc}(B)]$

$\subseteq \overline{\gamma}^{Inc}(A) \cap \overline{\gamma}^{Inc}(B)$.

(3) $\overline{\gamma}^{Inc}(A \cup B) = (A \cup B) \cup [\overline{R}^{Inc}(\underline{R}_{Inc}(A \cup B) \cup \underline{R}_{Inc}(\overline{R}^{Inc}(A \cup B)))]$

$\supseteq (A \cup B) \cup [\overline{R}^{Inc}(\underline{R}_{Inc}(A) \cup \underline{R}_{Inc}(B)) \cup \underline{R}_{Inc}(\overline{R}^{Inc}(A) \cup \overline{R}^{Inc}(B)))]$

$\supseteq (A \cup B) \cup [\overline{R}^{Inc}(\underline{R}_{Inc}(A)) \cup \overline{R}^{Inc}(\underline{R}_{Inc}(B)) \cup \underline{R}_{Inc}(\overline{R}^{Inc}(A)) \cup \underline{R}_{Inc}(\overline{R}^{Inc}(B))]$

$\supseteq A \cup [\overline{R}^{Inc}(\underline{R}_{Inc}(A)) \cup \underline{R}_{Inc}(\overline{R}^{Inc}(A))] \cup B \cup [\overline{R}^{Inc}(\underline{R}_{Inc}(B)) \cup \underline{R}_{Inc}(\overline{R}^{Inc}(B))]$

$$\supseteq \overline{\gamma}^{Inc}(A) \cup \overline{\gamma}^{Inc}(B).$$

One can prove the case between parentheses.

**Proposition 3.3.** Let $(U, \tau_R, \rho)$ be a GOTAS and $A, B \subseteq U$. Then

(1) $A \subseteq B \rightarrow \underline{\gamma}_{Inc}(A) \subseteq \underline{\gamma}_{Inc}(B)$ ( $A \subseteq B \rightarrow \underline{\gamma}_{Dec}(A) \subseteq \underline{\gamma}_{Dec}(B)$ ).

(2) $\underline{\gamma}_{Inc}(A \cap B) \subseteq \underline{\gamma}_{Inc}(A) \cap \underline{\gamma}_{Inc}(B)$ ($\underline{\gamma}_{Dec}(A \cap B) \subseteq \underline{\gamma}_{Dec}(A) \cap \underline{\gamma}_{Dec}(B)$).

(3) $\underline{\gamma}_{Inc}(A \cup B) \supseteq \underline{\gamma}_{Inc}(A) \cup \underline{\gamma}_{Inc}(B)$ ( $\underline{\gamma}_{Dec}(A \cup B) \supseteq \underline{\gamma}_{Dec}(A) \cup \underline{\gamma}_{Dec}(B)$ ).

**Proof.**

(1) Easy.

(2) $\underline{\gamma}_{Inc}(A \cap B) = (A \cap B) \cap [\overline{R}^{Inc}(\underline{R}_{Inc}(A \cap B) \cup \underline{R}_{Inc}(\overline{R}^{Inc}(A \cap B)))]$

$$\subseteq (A \cap B) \cap [\overline{R}^{Inc}(\underline{R}_{Inc}(A) \cap \underline{R}_{Inc}(B)) \cup \underline{R}_{Inc}(\overline{R}^{Inc}(A) \cap \overline{R}^{Inc}(B)))]$$

$$\subseteq (A \cap B) \cap [\overline{R}^{Inc}(\underline{R}_{Inc}(A) \cap \overline{R}^{Inc} \underline{R}_{Inc}(B)) \cup \underline{R}_{Inc}(\overline{R}^{Inc}(A) \cap \underline{R}_{Inc} \overline{R}^{Inc}(B)))]$$

$$\subseteq A \cap [\overline{R}^{Inc}(\underline{R}_{Inc}(A) \cup \underline{R}_{Inc}(\overline{R}^{Inc}(A))] \cap B \cap [\overline{R}^{Inc}(\underline{R}_{Inc}(B) \cup \underline{R}_{Inc}(\overline{R}^{Inc}(B))]$$

$$\subseteq \underline{\gamma}_{Inc}(A) \cap \underline{\gamma}_{Inc}(B).$$

(3) $\underline{\gamma}_{Inc}(A \cup B) = (A \cup B) \cap [\overline{R}^{Inc}(\underline{R}_{Inc}(A \cup B) \cup \underline{R}_{Inc}(\overline{R}^{Inc}(A \cup B)))]$

$$\supseteq (A \cup B) \cap [\overline{R}^{Inc}(\underline{R}_{Inc}(A) \cup \underline{R}_{Inc}(B)) \cup \underline{R}_{Inc}(\overline{R}^{Inc}(A) \cup \overline{R}^{Inc}(B)))]$$

$$\supseteq (A \cup B) \cap [\overline{R}^{Inc}(\underline{R}_{Inc}(A)) \cup \overline{R}^{Inc} \underline{R}_{Inc}(B)) \cup \underline{R}_{Inc}(\overline{R}^{Inc}(A)) \cup \underline{R}_{Inc} \overline{R}^{Inc}(B)))]$$

$$\supseteq A \cap [\overline{R}^{Inc}(\underline{R}_{Inc}(A) \cup \underline{R}_{Inc}(\overline{R}^{Inc}(A))] \cup B \cap [\overline{R}^{Inc}(\underline{R}_{Inc}(B) \cup \underline{R}_{Inc}(\overline{R}^{Inc}(B))]$$

$$\supseteq \underline{\gamma}_{Inc}(A) \cup \underline{\gamma}_{Inc}(B).$$

One can prove the case between parentheses.

**Proposition 3.4.** Let $(U, \tau_R, \rho)$ be a GOTAS and $A, B \subseteq U$. If $A$ is $R$ – increasing (resp. decreasing) exact then $A$ is $R$ – increasing (resp. decreasing) $\gamma$ exact.

**Proof.**

Let $A$ be $R$ – increasing exact. Then $\overline{R}^{Inc}(A) = \underline{R}_{Inc}(A)$, thus $\overline{\gamma}^{Inc}(A) = \overline{R}^{Inc}(A)$ and $\underline{\gamma}_{Inc}(A) = \underline{R}_{Inc}(A)$. Therefore $\overline{\gamma}^{Inc}(A) = \underline{\gamma}_{Inc}(A)$.

One can prove the case between parentheses.

$R$ – increasing (resp. decreasing) exact $\longrightarrow$ $R$ – increasing (resp. decreasing) $\gamma$ exact

**Proposition 3.5.** Let $(U, \tau_R, \rho)$ be a GOTAS and $A \subseteq U$. Then
$\underline{R}_{Inc}(A) \subseteq \underline{\gamma}_{Inc}(A)$ ($\underline{R}_{Dec}(A) \subseteq \underline{\gamma}_{Dec}(A)$).

**Proof.** Since $\underline{R}_{Inc}(A) \subseteq A$ and $\underline{R}_{Inc}(A) \subseteq \overline{R}^{Inc}(\underline{R}_{Inc}(A))$, then $\underline{R}_{Inc}(A) \subseteq \overline{R}^{Inc}(\underline{R}_{Inc}(A)) \cup \underline{R}_{Inc}(\overline{R}^{Inc}(A))$. Therefore $\underline{R}_{Inc}(A) \subseteq A \cap [\overline{R}^{Inc}(\underline{R}_{Inc}(A)) \cup \underline{R}_{Inc}(\overline{R}^{Inc}(A))]$. Thus $\underline{R}_{Inc}(A) \subseteq \underline{\gamma}_{Inc}(A)$.

One can prove the case between parentheses.

**Proposition 3.6.** Let $(U, \tau_R, \rho)$ be a GOTAS and $A \subseteq U$. Then
$\overline{\gamma}^{Inc}(A) \subseteq \overline{R}^{Inc}(A)$ ($\overline{\gamma}^{Dec}(A) \subseteq \overline{R}^{Dec}(A)$).

**Proof.** Since $A \subseteq \overline{R}^{Inc}(A)$ and $\underline{R}_{Inc}(A) \subseteq A \subseteq \overline{R}^{Inc}(A)$, then $\overline{R}^{Inc}(\underline{R}_{Inc}(A)) \subseteq A \subseteq \overline{R}^{Inc}(A)$. Thus $\overline{R}^{Inc}(\underline{R}_{Inc}(A)) \cup \underline{R}_{Inc}(\overline{R}^{Inc}(A)) \subseteq \overline{R}^{Inc}(A)$. Therefore $A \cup \overline{R}^{Inc}(\underline{R}_{Inc}(A)) \cup \underline{R}_{Inc}(\overline{R}^{Inc}(A)) \subseteq \overline{R}^{Inc}(A)$. Hence $\overline{\gamma}^{Inc}(A) \subseteq \overline{R}^{Inc}(A)$.

**Proposition 3.7.** Let $(U, \tau_R, \rho)$ be a GOTAS and $A \subseteq U$. Then
$\underline{P}_{Inc}(A) \subseteq \underline{\gamma}_{Inc}(A)$ ($\underline{P}_{Dec}(A) \subseteq \underline{\gamma}_{Dec}(A)$).

**Proof.** Let $x \in \underline{P}_{Inc}(A) = A \cap \underline{R}_{Inc}(\overline{R}^{Inc}(A))$. Then $x \in A$ and $\underline{R}_{Inc}(\overline{R}^{Inc}(A))$. Therefore $x \in A$ and $[x \in \overline{R}^{Inc}(\underline{R}_{Inc}(A)) \cup x \in \underline{R}_{Inc}(\overline{R}^{Inc}(A))]$. Thus $x \in A \cap [\overline{R}^{Inc}(\underline{R}_{Inc}(A)) \cup \underline{R}_{Inc}(\overline{R}^{Inc}(A))] = \underline{\gamma}_{Inc}(A)$. Hence $\underline{P}_{Inc}(A) \subseteq \underline{\gamma}_{Inc}(A)$.

One can prove the case between parentheses.

**Proposition 3.8.** Let $(U, \tau_R, \rho)$ be a GOTAS and $A \subseteq U$. Then
$\underline{S}_{Inc}(A) \subseteq \underline{\gamma}_{Inc}(A)$ ($\underline{S}_{Dec}(A) \subseteq \underline{\gamma}_{Dec}(A)$).

**Proof.** Let $x \in \underline{S}_{Inc}(A) = A \cap \overline{R}^{Inc}(\underline{R}_{Inc}(A))$. Then $x \in A$ and $\overline{R}^{Inc}(\underline{R}_{Inc}(A))$. Therefore $x \in A$ and $[x \in \overline{R}^{Inc}(\underline{R}_{Inc}(A))$ or $x \in \underline{R}_{Inc}(\overline{R}^{Inc}(A))]$. Thus $x \in A \cap [\overline{R}^{Inc}(\underline{R}_{Inc}(A)) \cup \underline{R}_{Inc}(\overline{R}^{Inc}(A))] = \underline{\gamma}_{Inc}(A)$. Hence $\underline{S}_{Inc}(A) \subseteq \underline{\gamma}_{Inc}(A)$.

One can prove the case between parentheses.

**Proposition 3.9.** Let $(U, \tau_R, \rho)$ be a GOTAS and $A \subseteq U$. Then $\overline{P}^{Inc}(A) \subseteq \overline{\gamma}^{Inc}(A) (\overline{P}^{Dec}(A) \subseteq \overline{\gamma}^{Dec}(A))$.

**Proof.** Let $x \in \overline{P}^{Inc}(A) = A \cup \overline{R}^{Inc}(\underline{R}_{Inc}(A))$. Then $x \in A$ and $\overline{R}^{Inc}(\underline{R}_{Inc}(A))$. Therefore $x \in A \cup [\overline{R}^{Inc}(\underline{R}_{Inc}(A)) \cup \underline{R}_{Inc}(\overline{R}^{Inc}(A))]$. Thus $\overline{P}^{Inc}(A) \subseteq \overline{\gamma}^{Inc}(A)$.

**Proposition 3.10.** Let $(U, \tau_R, \rho)$ be a GOTAS and $A \subseteq U$. Then $\overline{\beta}^{Inc}(A) \subseteq \overline{P}^{Inc}(A) (\overline{\beta}^{Dec}(A) \subseteq \overline{P}^{Dec}(A))$.
Proof. Omitted.

**Definition 3.11.** Let $(U, \tau_R, \rho)$ be a GOTAS and $A \subseteq U$. We define:

(1) $\underline{\beta}_{Inc}(A) = A \cap \overline{R}^{Inc}(\underline{R}_{Inc}(\overline{R}^{Inc}(A)))$, $\underline{\beta}_{Inc}(A)$ is called $R$–increasing $\beta$ lower.

(2) $\overline{\beta}^{Inc}(A) = A \cup \underline{R}_{Inc}(\overline{R}^{Inc}(\underline{R}_{Inc}(A))$, $\overline{\beta}^{Inc}(A)$ is called $R$–increasing $\beta$ upper.

(3) $\underline{\beta}_{Dec}(A) = A \cap \overline{R}^{Dec}(\underline{R}_{Dec}(\overline{R}^{Dec}(A)))$, $\underline{\beta}_{Dec}(A)$ is called $R$–decreasing $\beta$ lower.

(4) $\overline{\beta}^{Dec}(A) = A \cup \underline{R}_{Dec}(\overline{R}^{Dec}(\underline{R}_{Dec}(A))$, $\overline{\beta}^{Dec}(A)$ is called $R$–decreasing $\beta$ upper.

$A$ is $R$–increasing (decreasing) $\beta$ exact if $\underline{\beta}_{Inc}(A) = \overline{\beta}^{Inc}(A)$ (resp. $\underline{\beta}_{Dec}(A) = \overline{\beta}^{Dec}(A)$), otherwise $A$ is $R$–increasing (decreasing) $\beta$ rough.

**Proposition 3.12.** Let $(U, \tau_R, \rho)$ be a GOTAS and $A, B \subseteq U$. Then

(1) $A \subseteq B \rightarrow \overline{\beta}^{Inc}(A) \subseteq \overline{\beta}^{Inc}(B) (A \subseteq B \rightarrow \overline{\beta}^{Dec}(A) \subseteq \overline{\beta}^{Dec}(B))$.

(2) $\overline{\beta}^{Inc}(A \cap B) \subseteq \overline{\beta}^{Inc}(A) \cap \overline{\beta}^{Inc}(B) (\overline{\beta}^{Dec}(A \cap B) \subseteq \overline{\beta}^{Dec}(A) \cap \overline{\beta}^{Dec}(B))$.

(3) $\overline{\beta}^{Inc}(A \cup B) \supseteq \overline{\beta}^{Inc}(A) \cup \overline{\beta}^{Inc}(B) (\overline{\beta}^{Dec}(A \cup B) \supseteq \overline{\beta}^{Dec}(A) \cup \overline{\beta}^{Dec}(B))$.

**Proof.**

(1) Omitted.

(2) $\overline{\beta}^{Inc}(A \cap B) = (A \cap B) \cup \underline{R}_{Inc}(\overline{R}^{Inc}(\underline{R}_{Inc}(A \cap B)))$

$$= (A \cap B) \cup \underline{R}_{Inc}(\overline{R}^{Inc}(\underline{R}_{Inc}(A) \cap \underline{R}_{Inc}(B)))$$

$$\subseteq (A \cap B) \cup \underline{R}_{Inc}(\overline{R}^{Inc}(\underline{R}_{Inc}(A) \cap \overline{R}^{Inc}(\underline{R}_{Inc}(B))))$$

$$\subseteq (A \cap B) \cup \underline{R}_{Inc}(\overline{R}^{Inc}(\underline{R}_{Inc}(A)) \cap \underline{R}_{Inc}(\overline{R}^{Inc}(\underline{R}_{Inc}(B))))$$

$$\subseteq A \cup \underline{R}_{Inc}(\overline{R}^{Inc}(\underline{R}_{Inc}(A))) \cap B \cup \underline{R}_{Inc}(\overline{R}^{Inc}(\underline{R}_{Inc}(B)))$$

$$\subseteq \overline{\beta}^{Inc}(A) \cap \overline{\beta}^{Inc}(B).$$

(3) $\overline{\beta}^{Inc}(A \cup B) = (A \cup B) \cup \underline{R}_{Inc}(\overline{R}^{Inc}(\underline{R}_{Inc}(A \cup B)))$

$$= (A \cup B) \cup \underline{R}_{Inc}(\overline{R}^{Inc}(\underline{R}_{Inc}(A) \cup \underline{R}_{Inc}(B))$$

$$\supseteq (A \cup B) \cup \underline{R}_{Inc}(\overline{R}^{Inc}(\underline{R}_{Inc}(A) \cup \overline{R}^{Inc}(\underline{R}_{Inc}(B)))$$

$$\supseteq (A \cup B) \cup \underline{R}_{Inc}(\overline{R}^{Inc}(\underline{R}_{Inc}(A) \cup \underline{R}_{Inc}(\overline{R}^{Inc}(\underline{R}_{Inc}(B))))$$

$$\supseteq (A \cup \underline{R}_{Inc}(\overline{R}^{Inc}(\underline{R}_{Inc}(A))) \cup (B \cup \underline{R}_{Inc}(\overline{R}^{Inc}(\underline{R}_{Inc}(B))))$$

$$\supseteq \overline{\beta}^{Inc}(A) \cup \overline{\beta}^{Inc}(B).$$

One can prove the case between parentheses.

**Proposition 3.13.** Let $(U, \tau_R, \rho)$ be a GOTAS and $A, B \subseteq U$. Then
  (1) $A \subseteq B \to \underline{\beta}_{Inc}(A) \subseteq \underline{\beta}_{Inc}(B)$ ( $A \subseteq B \to \underline{\beta}_{Dec}(A) \subseteq \underline{\beta}_{Dec}(B)$ ).
  (2) $\underline{\beta}_{Inc}(A \cap B) \subseteq \underline{\beta}_{Inc}(A) \cap \underline{\beta}_{Inc}(B)$ ( $\underline{\beta}_{Dec}(A \cap B) \subseteq \underline{\beta}_{Dec}(A) \cap \underline{\beta}_{Dec}(B)$ ).
  (3) $\underline{\beta}_{Inc}(A \cup B) \supseteq \underline{\beta}_{Inc}(A) \cup \underline{\beta}_{Inc}(B)$ ( $\underline{\gamma}_{Dec}(A \cup B) \supseteq \underline{\beta}_{Dec}(A) \cup \underline{\beta}_{Dec}(B)$ ).

**Proof.**

(1) Easy.

(2) $\underline{\beta}_{Inc}(A \cap B) = (A \cap B) \cap \overline{R}^{Inc}(\underline{R}_{Inc}(\overline{R}^{Inc}(A \cap B)))$

$$\subseteq (A \cap B) \cap \overline{R}^{Inc}(\underline{R}_{Inc}(\overline{R}^{Inc}(A) \cap \overline{R}^{Inc}(B)))$$

$$\subseteq (A \cap B) \cap \overline{R}^{Inc}(\underline{R}_{Inc}(\overline{R}^{Inc}(A) \cap \overline{R}^{Inc}(\underline{R}_{Inc}(\overline{R}^{Inc}(B)))$$

$$\subseteq A \cap \overline{R}^{Inc}(\underline{R}_{Inc}(\overline{R}^{Inc}(A) \cap B \cap \overline{R}^{Inc}(\underline{R}_{Inc}(\overline{R}^{Inc}(B)))$$

$$\subseteq \underline{\beta}_{Inc}(A) \cap \underline{\beta}_{Inc}(B).$$

(3) $\underline{\beta}_{Inc}(A \cup B) = (A \cup B) \cap \overline{R}^{Inc}(\underline{R}_{Inc}(\overline{R}^{Inc}(A \cup B)))$

$= (A \cup B) \cap \overline{R}^{Inc}(\underline{R}_{Inc}(\overline{R}^{Inc}(A \cup B)))$

$\supseteq (A \cup B) \cap \overline{R}^{Inc}(\underline{R}_{Inc}(\overline{R}^{Inc}(A) \cup \overline{R}^{Inc}(B)))$

$\supseteq (A \cup B) \cap \overline{R}^{Inc}(\underline{R}_{Inc}(\overline{R}^{Inc}(A) \cup \overline{R}^{Inc}(\underline{R}_{Inc}(\overline{R}^{Inc}(B))))$

$\supseteq A \cap \overline{R}^{Inc}(\underline{R}_{Inc}(\overline{R}^{Inc}(A) \cup B \cap \overline{R}^{Inc}(\underline{R}_{Inc}(\overline{R}^{Inc}(B))))$

$\supseteq \underline{\beta}_{Inc}(A) \cup \underline{\beta}_{Inc}(B)$.

One can prove the case between parentheses.

**Proposition 3.14.** Let $(U, \tau_R, \rho)$ be a GOTAS and $A, B \subseteq U$. If $A$ is $R-$ increasing (resp. decreasing) exact then $A$ is $\beta-$ increasing (resp. decreasing) exact.

**Proof.**

Let $A$ be $R-$ increasing exact. Then $\overline{R}^{Inc}(A) = \underline{R}_{Inc}(A)$. Therefore $\overline{\beta}^{Inc}(A) = \overline{R}^{Inc}(A)$, $\underline{\beta}_{Inc}(A) = \underline{R}_{Inc}(A)$. Thus $\overline{\beta}^{Inc}(A) = \underline{\beta}_{Inc}(A)$. Hence $A$ is $R-$ increasing $\beta$ exact.

One can prove the case between parentheses.

$R-$ increasing (resp. decreasing) exact $\longrightarrow$ $R-$ increasing (resp. decreasing) $\beta$ exact

**Proposition 3.15.** Let $(U, \tau_R, \rho)$ be a GOTAS and $A \subseteq U$. Then
$\underline{R}_{Inc}(A) \subseteq \underline{\beta}_{Inc}(A)\,(\underline{R}_{Dec}(A) \subseteq \underline{\beta}_{Dec}(A))$

**Proof.** Since $\underline{R}_{Inc}(A) \subseteq A \subseteq \overline{R}^{Inc}(A)$ and $\underline{R}_{Inc}(A) \subseteq \underline{R}_{Inc}(\overline{R}^{Inc}(A))$. Then
$\underline{R}_{Inc}(A) \subseteq \overline{R}^{Inc}(\underline{R}_{Inc}(A)) \subseteq \overline{R}^{Inc}(\underline{R}_{Inc}(\overline{R}^{Inc}(A)))$. Therefore
$\underline{R}_{Inc}(A) \subseteq A \cap [\overline{R}^{Inc}(\underline{R}_{Inc}(\overline{R}^{Inc}(A)))]$. Thus $\underline{R}_{Inc}(A) \subseteq \underline{\beta}_{Inc}(A)$.

One can prove the case between parentheses.

**Proposition 3.16.** Let $(U, \tau_R, \rho)$ be a GOTAS and $A \subseteq U$. Then
$\overline{\beta}^{Inc}(A) \subseteq \overline{R}^{Inc}(A)\,(\overline{\beta}^{Dec}(A) \subseteq \overline{R}^{Dec}(A))$.

**Proof.** Since $A \subseteq \overline{R}^{Inc}(A)$ and $\underline{R}_{Inc}(A) \subseteq \overline{R}^{Inc}(A)$. Then $\overline{R}^{Inc}(\underline{R}_{Inc}(A)) \subseteq \overline{R}^{Inc}(A)$. Thus $\underline{R}_{Inc}(\overline{R}^{Inc}(\underline{R}_{Inc}(A)) \subseteq \underline{R}_{Inc}(\overline{R}^{Inc}(A)) \subseteq \overline{R}^{Inc}(A)$. Therefore $A \cup \underline{R}_{Inc}(\overline{R}^{Inc}(\underline{R}_{Inc}(A)) \subseteq \overline{R}^{Inc}(A)$. Hence $\overline{\beta}^{Inc}(A) \subseteq \overline{R}^{Inc}(A)$.

**Definition 3.17.** Let $(U, \tau_R, \rho)$ be a GOTAS and $A \subseteq U$. Then

(1) $B_{jInc}(A) = \overline{j}^{Inc}(A) - \underline{j}_{Inc}(A)$ (resp. $B_{jDec}(A) = \overline{j}^{Dec}(A) - \underline{j}_{Dec}(A)$),

is increasing (resp. decreasing) $j$ boundary region.

(2) $Pos_{jInc}(A) = \underline{j}_{Inc}(A)$ (resp. $Pos_{jDec}(A) = \underline{j}_{Dec}(A)$),

is increasing (resp. decreasing) $j$ positive region.

(3) $Neg_{jInc}(A) = U - \overline{j}^{Dec}(A)$ (resp. $Neg_{jDec}(A) = U - \overline{j}^{Inc}(A)$),

is increasing (resp. decreasing) $j$ negative region. Where $\underline{j}_{Inc}$ the near lower approximations s.t. $j \in \{\beta, \gamma\}$.

**Proposition 3.18.** Let $(U, \tau_R, \rho)$ be a GOTAS and $A, B \subseteq U$. Then

(1) $Neg_{\gamma Inc}(A \cup B) \subseteq Neg_{\gamma Inc}(A) \cup Neg_{\gamma Inc}(B)$

($Neg_{\gamma Dec}(A \cup B) \subseteq Neg_{\gamma Dec}(A) \cup Neg_{\gamma Dec}(B)$).

(2) $Neg_{\gamma Inc}(A \cap B) \supseteq Neg_{\gamma Inc}(A) \cap Neg_{\gamma Inc}(B)$

($Neg_{\gamma Dec}(A \cap B) \supseteq Neg_{\gamma Dec}(A) \cap Neg_{\gamma Dec}(B)$).

**Proof.**

(1) $Neg_{\gamma Inc}(A \cup B) = U - [(A \cup B) \cup [\overline{R}^{Dec} \underline{R}_{Dec}(A \cup B) \cup \underline{R}_{Dec}\overline{R}^{Dec}(A \cup B)]]$

$\subseteq U - [(A \cup B) \cup [\overline{R}^{Dec}(\underline{R}_{Dec}(A) \cup \underline{R}_{Dec}(B)) \cup \underline{R}_{Dec}(\overline{R}^{Dec}(A) \cup \overline{R}^{Dec}(B))]]$

$\subseteq U - [(A \cup B) \cup [\overline{R}^{Dec}\underline{R}_{Dec}(A) \cup \overline{R}^{Dec}\underline{R}_{Dec}(B) \cup \underline{R}_{Dec}\overline{R}^{Dec}(A) \cup \underline{R}_{Dec}\overline{R}^{Dec}(B)]]$

$\subseteq U - [A \cup [\overline{R}^{Dec}\underline{R}_{Dec}(A) \cup \underline{R}_{Dec}\overline{R}^{Dec}(A)] \cup [B \cup [\overline{R}^{Dec}\underline{R}_{Dec}(B) \cup \underline{R}_{Dec}\overline{R}^{Dec}(B)]]$

$\subseteq U - [A \cup [\overline{R}^{Dec}\underline{R}_{Dec}(A) \cup \underline{R}_{Dec}\overline{R}^{Dec}(A)] \cap U - B \cup [\overline{R}^{Dec}\underline{R}_{Dec}(B) \cup \underline{R}_{Dec}\overline{R}^{Dec}(B)]$

$\subseteq Neg_{\gamma Inc}(A) \cap Neg_{\gamma Inc}(B)$.

(2) $Neg_{\gamma Inc}(A \cap B) = U - [(A \cap B) \cup [\overline{R}^{Dec}\underline{R}_{Dec}(A \cap B) \cup \underline{R}_{Dec}\overline{R}^{Dec}(A \cap B)]]$

$$\supseteq U - [(A \cap B) \cup [\overline{R}^{Dec}(\underline{R}_{Dec}(A) \cap \underline{R}_{Dec}(B)) \cup \underline{R}_{Dec}(\overline{R}^{Dec}(A) \cap \overline{R}^{Dec}(B))]]$$

$$\supseteq U - [(A \cap B) \cup [\overline{R}^{Dec}\underline{R}_{Dec}(A) \cap \overline{R}^{Dec}\underline{R}_{Dec}(B) \cup \underline{R}_{Dec}\overline{R}^{Dec}(A) \cap \underline{R}_{Dec}\overline{R}^{Dec}(B)]]$$

$$\supseteq U - [A \cup [\overline{R}^{Dec}\underline{R}_{Dec}(A) \cup \underline{R}_{Dec}\overline{R}^{Dec}(A)] \cap [B \cup [\overline{R}^{Dec}\underline{R}_{Dec}(B) \cup \underline{R}_{Dec}\overline{R}^{Dec}(B)]]$$

$$\supseteq U - [A \cup [\overline{R}^{Dec}\underline{R}_{Dec}(A) \cup \underline{R}_{Dec}\overline{R}^{Dec}(A)] \cup U - [B \cup [\overline{R}^{Dec}\underline{R}_{Dec}(B) \cup \underline{R}_{Dec}\overline{R}^{Dec}(B)]]$$

$$\supseteq Neg_{\gamma Inc}(A) \cup Neg_{\gamma Inc}(B).$$

One can prove the case between parentheses.

**Proposition 3.19.** Let $(U, \tau_R, \rho)$ be a GOTAS and $A, B \subseteq U$. Then

(1) $Neg_{\beta Inc}(A \cup B) \subseteq Neg_{\beta Inc}(A) \cup Neg_{\beta Inc}(B)$

$(Neg_{\beta Dec}(A \cup B) \subseteq Neg_{\beta Dec}(A) \cup Neg_{\beta Dec}(B)).$

(2) $Neg_{\beta Inc}(A \cap B) \supseteq Neg_{\beta Inc}(A) \cap Neg_{\beta Inc}(B)$

$(Neg_{\beta Dec}(A \cap B) \supseteq Neg_{\beta Dec}(A) \cap Neg_{\beta Dec}(B)).$

**Proof.**

(1) $Neg_{\beta Inc}(A \cup B) = U - [(A \cup B) \cup \underline{R}_{Dec}(\overline{R}^{Dec}\underline{R}_{Dec}(A \cup B)]$

$$\subseteq U - [(A \cup B) \cup \underline{R}_{Dec}(\overline{R}^{Dec}(\underline{R}_{Dec}(A) \cup \underline{R}_{Dec}(B)))]$$

$$\subseteq U - [(A \cup B) \cup \underline{R}_{Dec}(\overline{R}^{Dec}(\underline{R}_{Dec}(A) \cup \overline{R}^{Dec}(\underline{R}_{Dec}(B))))]$$

$$\subseteq U - [(A \cup B) \cup (\underline{R}_{Dec}(\overline{R}^{Dec}(\underline{R}_{Dec}(A) \cup \underline{R}_{Dec}\overline{R}^{Dec}(\underline{R}_{Dec}(B))))]$$

$$\subseteq U - [A \cup \underline{R}_{Dec}(\overline{R}^{Dec}(\underline{R}_{Dec}(A))) \cup (B \cup \underline{R}_{Dec}\overline{R}^{Dec}(\underline{R}_{Dec}(B)))]$$

$$\subseteq Neg_{\beta Inc}(A) \cap Neg_{\beta Inc}(B).$$

(2) $Neg_{\beta Inc}(A \cap B) = U - [(A \cap B) \cup \underline{R}_{Dec}(\overline{R}^{Dec}\underline{R}_{Dec}(A \cap B))]$

$$= U - [(A \cap B) \cup \underline{R}_{Dec}(\overline{R}^{Dec}(\underline{R}_{Dec}(A) \cap \underline{R}_{Dec}(B)))]$$

$$\supseteq U - [(A \cap B) \cup \underline{R}_{Dec}(\overline{R}^{Dec}(\underline{R}_{Dec}(A) \cap \underline{R}_{Dec}(\overline{R}^{Dec}(\underline{R}_{Dec}(B))))]$$

$$\supseteq U - [A \cup \underline{R}_{Dec}(\overline{R}^{Dec}(\underline{R}_{Dec}(A))) \cap B \cup \underline{R}_{Dec}(\overline{R}^{Dec}(\underline{R}_{Dec}(B)))]$$

$$\supseteq U - A \cup \underline{R}_{Dec}(\overline{R}^{Dec}(\underline{R}_{Dec}(A)) \cup U - B \cup \underline{R}_{Dec}(\overline{R}^{Dec}(\underline{R}_{Dec}(B)))$$

$$\supseteq Neg_{\beta Inc}(A) \cup Neg_{\beta Inc}(B).$$

One can prove the case between parentheses.

**Proposition 3.20.** Let $(U, \tau_R, \rho)$ be a GOTAS and $A \subseteq U$. Then
$$\underline{S}_{Inc}(A) \subseteq \underline{\gamma}_{Inc}(A) \subseteq \underline{\beta}_{Inc}(A) \, (\underline{S}_{Dec}(A) \subseteq \underline{\gamma}_{Dec}(A) \subseteq \underline{\beta}_{Dec}(A)).$$

**Proof.**

Let $x \in \underline{S}_{Inc}(A)$. Then $x \in \overline{R}^{Inc}(\underline{R}_{Inc}(A))$. Therefore $x \in \overline{R}^{Inc}(\underline{R}_{Inc}(A)) \cup \underline{R}_{Inc}(\overline{R}^{Inc}(A))$.

Thus $x \in A \cap [\overline{R}^{Inc}(\underline{R}_{Inc}(A)) \cup \underline{R}_{Inc}(\overline{R}^{Inc}(A))]$ and thus $x \in \underline{\gamma}_{Inc}(A)$.

Hence $\underline{S}_{Inc}(A) \subseteq \underline{\gamma}_{Inc}(A)$ \hspace{2em} (1).

Since $x \in \underline{R}_{Inc}(A)$, then $x \in \overline{R}^{Inc}(A)$. Therefore $x \in \underline{R}_{Inc}(\overline{R}^{Inc}(A))$.

Thus $x \in \overline{R}^{Inc}(\underline{R}_{Inc}(\overline{R}^{Inc}(A)))$, and thus $x \in A \cap \overline{R}^{Inc}(\underline{R}_{Inc}(\overline{R}^{Inc}(A)))$. Hence
$$x \in \underline{\beta}_{Inc}(A) \hspace{2em} (2)$$

From (1) and (2) we have,
$$\underline{S}_{Inc}(A) \subseteq \underline{\gamma}_{Inc}(A) \subseteq \underline{\beta}_{Inc}(A).$$

One can prove the case between parentheses.

**Proposition 3.21.** Let $(U, \tau_R, \rho)$ be a GOTAS and $A \subseteq U$. Then
$$\overline{\beta}^{Inc}(A) \subseteq \overline{\gamma}^{Inc}(A) \subseteq \overline{S}^{Inc}(A) \, (\overline{\beta}^{Dec}(A) \subseteq \overline{\gamma}^{Dec}(A) \subseteq \overline{S}^{Dec}(A)).$$

**Proof.**

Let $x \in \overline{\beta}^{Inc}(A)$. Then $x \in A \cup \underline{R}_{Inc}(\overline{R}^{Inc}(\underline{R}_{Inc}(A)))$. Therefore

$x \in A$ or $x \in \underline{R}_{Inc}(\overline{R}^{Inc}(\underline{R}_{Inc}(A)))$. Thus $x \in A$ or $x \in \overline{R}^{Inc}(\underline{R}_{Inc}(A))$. So $x \in A \cup \overline{R}^{Inc}(\underline{R}_{Inc}(A))$, and so $x \in A \cup [\overline{R}^{Inc}(\underline{R}_{Inc}(A)) \cup \underline{R}_{Inc}(\overline{R}^{Inc}(A))]$. Thus $x \in \overline{\gamma}^{Inc}(A)$. Hence $\overline{\beta}^{Inc}(A) \subseteq \overline{\gamma}^{Inc}(A)$. \hspace{2em} (1)

Since $x \in \overline{\gamma}^{Inc}(A)$, $x \in A$ or $x \in \underline{R}_{Inc}(\overline{R}^{Inc}(A))$, then $x \in A \cup \underline{R}_{Inc}(\overline{R}^{Inc}(A))$. Therefore
$$x \in \overline{S}^{Inc}(A) \hspace{2em} (2)$$

From (1) and (2) we have, $\overline{\beta}^{Inc}(A) \subseteq \overline{\gamma}^{Inc}(A) \subseteq \overline{S}^{Inc}(A)$.

One can prove the case between parentheses.

**Definition 3.22.** Let $(U, \tau_R, \rho)$ be a GOTAS and $A$ is a non-empty finite subset of $U$. Then the increasing (decreasing) j accuracy of a finite non-empty subset $A$ of $U$ is given by:

$$\eta_{jInc}(A) = \frac{\left|\underline{j}_{Inc}(A)\right|}{\left|\overline{j}^{Inc}(A)\right|}, \quad j \in \{\beta, \gamma\}.$$

**Proposition 3.23.** Let $(U, \tau_R, \rho)$ be a GOTAS and $A$ non-empty finite subset of $U$. Then we have $\eta(A) \leq \eta_{jInc}(A)$ $(\eta(A) \leq \eta_{jDec}(A))$, for all $j \in \{\beta, \gamma\}$, where $\eta(A) = \frac{|\underline{R}(A)|}{|\overline{R}(A)|}$.

**Proof.** Omitted.

**Example 3.24.** Let $U = \{a,b,c,d\}$, $U/R = \{\{a\},\{a,b\},\{c,d\}\}$, $\tau_R = \{U, \phi, \{a,b\}, \{c,d\}, \{a\}, \{a,d,c\}\}$, $\tau_R^C = \{U, \phi, \{c,d\}, \{a,b\}, \{b,c,d\}, \{b\}\}$ and $\rho = \{(a,a),(b,b),(c,c),(d,d),(a,b),(b,d),(a,d),(a,c),(c,d)\}$.

For $A = \{a,c\}$, we have:

$\underline{R}_{Dec}(A) = \{a\}$, $\overline{R}^{Dec}(\underline{R}_{Dec}(A)) = \{a,b\}$, $\overline{R}^{Dec}(A) = U$, $\underline{R}_{Dec}(\overline{R}^{Dec}(A)) = U$.

$\underline{S}_{Dec}(A) = \{a\}$, $\overline{S}^{Dec}(A) = U$, $B_{SDec}(A) = \{b,c,d\}$, $Neg_{SInc} = \phi$.

- - - - - - - - - - - - - - - - - -

$\underline{\gamma}_{Dec}(A) = A \cap U = A$, $\overline{\gamma}^{Dec}(A) = A \cup U = U$, $B_{\gamma Dec}(A) = \{b,d\}$, $Neg_{\gamma Inc} = \phi$.

- - - - - - - - - - - - - - - - - -

$\underline{\beta}_{Dec}(A) = A \cap U = A$, $\overline{\beta}^{Dec}(A) = \{a,b,c\}$, $B_{\beta Dec}(A) = \{b\}$, $Neg_{\beta Inc} = \{d\}$

**Proposition 3.25.** Let $(U, \tau_R, \rho)$ be a GOTAS and $A \subseteq U$. Then we have
$B_{\beta Inc}(A) \subseteq B_{\gamma Inc}(A) \subseteq B_{SInc}(A)$ $(B_{\beta Dec}(A) \subseteq B_{\gamma Dec}(A) \subseteq B_{SDec}(A))$
**Proof.** Omitted.

**Remark 3.26.** $B_{\gamma Inc}(A) \subseteq B_{RInc}(A) (B_{\gamma Dec}(A) \subseteq B_{RDec}(A))$.

**Remark 3.27.** $B_{\beta Inc}(A) \subseteq B_{RInc}(A) (B_{\beta Dec}(A) \subseteq B_{RDec}(A))$.

**Proposition 3.28.** Let $(U, \tau_R, \rho)$ be a GOTAS and $A$ be a non-empty finite subset of $U$. Then $\eta_{Inc}(A) \leq \eta_{\gamma Inc}(A) \leq \eta_{\beta Inc}(A)$ $(\eta_{Dec}(A) \leq \eta_{\gamma Dec}(A) \leq \eta_{\beta Dec}(A))$.
**Proof.** Omitted.

**Proposition 3.28.** Let $(U, \tau_R, \rho)$ be a GOTAS and $A \subseteq U$. Then $\underline{\gamma}_{Inc}(A) \subseteq \underline{\beta}_{Inc}(A) (\underline{\gamma}_{Dec}(A) \subseteq \underline{\beta}_{Dec}(A))$

**Proof.** Let $x \in \underline{\gamma}_{Inc}(A) = A \cap [\overline{R}^{Inc}(\underline{R}_{Inc}(A)) \cup \underline{R}_{Inc}(\overline{R}^{Inc}(A))]$. Then $x \in A$ and $x \in \overline{R}^{Inc}(\underline{R}_{Inc}(A)) \cup \underline{R}_{Inc}(\overline{R}^{Inc}(A))$. Therefore $x \in A$ and $[x \in \overline{R}^{Inc}(\underline{R}_{Inc}(A))$ or $x \in \underline{R}_{Inc}(\overline{R}^{Inc}(A))]$. Thus $x \in A$ and $x \in \underline{R}_{Inc}(\overline{R}^{Inc}(A))$ and thus $x \in A$ and $x \in \overline{R}^{Inc}(\underline{R}_{Inc}(\overline{R}^{Inc}(A)))$. Hence $x \in A \cap \overline{R}^{Inc}(\underline{R}_{Inc}(\overline{R}^{Inc}(A)))$. Therefore $\underline{\gamma}_{Inc}(A) \subseteq \underline{\beta}_{Inc}(A)$.

One can prove the case between parentheses.

## 4. Conclusion

Our results in this paper, became the results about of $\gamma$, $\beta$ approximation in [2] in case of $\rho$ is the equal relation. Also, The new approximation which we give became as Pawlak's approximation in case of $\rho$ is the equal relation and $R$ is the equivalence relation.